





\magnification\magstep1
\vsize=24.0truecm
\voffset-0.30truecm
\baselineskip13pt
\def\hatt{\widehat}
\def\dell{\partial}
\def\tilda{\widetilde}
\def\eps{\varepsilon}
\def\half{\hbox{$1\over2$}}

\font\csc=cmcsc10

\def\cstok#1{\leavevmode\thinspace\hbox{\vrule\vtop{\vbox{\hrule\kern1pt
	 \hbox{\vphantom{\tt/}\thinspace{\tt#1}\thinspace}}
	 \kern1pt\hrule}\vrule}\thinspace} 
\def\firkant{\cstok{\phantom{$\cdot$}}} 

\centerline{\bf On the last time and the number of times}
\centerline{\bf an estimator is more than $\eps$ from its target value} 

\medskip 
\centerline{\bf Nils Lid Hjort$^{a,b}$ and Grete Fenstad$^a$} 

\medskip
\centerline{\sl University of Oslo$^a$ 
		and Norwegian Computing Centre$^b$}

\smallskip
\centerline{\sl Statistical Research Report, April 1991}

\smallskip 
\centerline{\sl Department of Mathematics, University of Oslo}


\smallskip 

{{ \smallskip\narrower\noindent\baselineskip12pt 
{\csc Abstract}. 
Suppose $\hatt\theta_n$ is a strongly consistent 
estimator for $\theta_0$ in some i.i.d.~situation. 
Let $N_\eps$ and $Q_\eps$ be respectively 
the last $n$ and the total number of $n$ for which
$\hatt\theta_n$ is at least $\eps$ away from $\theta_0$. 
The limit distributions for $\eps^2N_\eps$ and $\eps^2Q_\eps$ 
as $\eps$ goes to zero are obtained under natural and weak conditions. 
The theory covers both 
parametric and nonparametric cases,
multi-dimensional parameters, and 
general distance functions.  
Our results are of probabilistic interest, and,
on the statistical side, suggest ways in which 
competing estimators can be compared. 
In particular several new optimality 
properties for the maximum likelihood estimator sequence
in parametric families are established. 
Another use of our results is ways of constructing 
sequential fixed-volume or shrinking-volume confidence sets,
as well as sequential tests with power 1. 
The paper also includes limit distribution results 
for the last $n$ and the number of $n$ 
for which the supremum distance $\|F_n-F\|\ge\eps$, 
where $F_n$ is the empirical distribution function. 
Other results are reached for
$\eps^{5/2}N_\eps$ and $\eps^{5/2}Q_\eps$ 
in the context of nonparametric density estimation, 
referring to the last time and the number of times 
where $|f_n(x)-f(x)|\ge\eps$. 
Finally it is shown that our results extend to several non-i.i.d.~situations.

\smallskip\noindent 
{\csc Key words:} 
\sl asymptotic optimality; 
asymptotic relative efficiency; 
density estimation; 
Glivenko--Cantelli;
the number of $\eps$-misses;
sequential confidence regions;
tests with power 1;
the last $n$
\smallskip}} 

\bigskip
{\bf 1. Introduction and summary.} 
Let $X_1,X_2,\ldots$ be a sequence of independent identically distributed 
(i.i.d.)~variables, and suppose $\hatt\theta_n$ is an estimator based
on the first $n$ observations which is strongly consistent
for some parameter $\theta_0$ of interest, 
i.e.~$\hatt\theta_n$ converges almost surely (a.s.) to $\theta_0$.
How large must $n$ be in order for 
$\hatt\theta_n$ to be very close to $\theta_0$? 

This natural question can be made precise in several different ways. 
(i) We can ask for an $m$ such that 
$${\rm Pr}\{|\hatt\theta_n-\theta_0|\le\eps\}\ge0.95 
	 \quad {\rm for\ all\ } n\ge m. \eqno(1.1)$$
An approximative answer to this question is readily given 
in the traditional cases where one has convergence in distribution of 
$\sqrt{n}(\hatt\theta_n-\theta_0)$ to some appropriate $N(0,\sigma_0^2)$.
Then $\sqrt{n}\eps/\sigma_0\ge1.96$ suffices, and we find 
$m\doteq1.96^2\sigma_0^2/\eps^2$ (and the assumption of 
strong consistency is not needed). 
(ii) We might ask for simultaneous 
closeness for all large $n$, with high enough probability, i.e.
$${\rm Pr}\{|\hatt\theta_n-\theta_0|\le\eps
	 {\rm\ for\ all\ }n\ge m\}
	 ={\rm Pr}\{\sup_{n\ge m}|\hatt\theta_n-\theta_0|\le\eps\}
		 \doteq 0.95, \eqno(1.2)$$
which also can be thought of as a requirement for a sequential
fixed-width confidence interval procedure. 
There is a finite $m$ solving this problem since 
$\sup_{n\ge m}|\hatt\theta_n-\theta_0|\rightarrow_p0$ 
when $\hatt\theta_n\rightarrow\theta_0$ a.s.
%
(iii) Thirdly we could study the random variable 
$$N_\eps=\sup\{n\ge1\colon|\hatt\theta_n-\theta_0|\ge\eps\}, \eqno(1.3)$$
which by the assumption of strong consistency is finite
with probability one. One has 
$${\rm Pr}\{\eps^2N_\eps\ge y\}
	 ={\rm Pr}\{N_\eps\ge m\}
	 ={\rm Pr}\{\sqrt{m}\sup_{n\ge m}|\hatt\theta_n
			 -\theta_0|\ge\sqrt{y_0}\}, \eqno(1.4)$$
in which $m=\langle y/\eps^2\rangle$
is the smallest integer $\ge y/\eps^2$, 
and $y_0=m\eps^2$ is close to $y$; $y_0-\eps^2<y\le y_0$. 
This shows that 
problems (ii) and (iii) are closely related; 
$\eps^2N_\eps$ has a limiting distribution if 
$\sqrt{m}\sup_{n\ge m}|\hatt\theta_n-\theta_0|$ has. 

While problem (i) is well-studied and solved problems (ii) and 
(iii) are virtually un-studied, presumably because 
the random variables they are concerned with depend upon
the full sequence of estimators and as such cannot be observed. 
They have nevertheless some immediate probabilistic and 
statistical appeal, and provide information 
about the speed of convergence of $\hatt\theta_n$ to its
target value. Stating or proving almost sure convergence  
touches basic chords in both probabilists and statisticians.
Since $\hatt\theta_n\rightarrow\theta_0$ a.s.~means 
nothing but that $N_\eps$ is a.s.~finite, 
it appears natural to inquire about its approximate size, 
for example via its approximate distribution 
and approximate expected value. 
Even Serfling's physician (1980, p.~49) is interested in (ii) and (iii). 
The `last $n$'  viewpoint also invites two competing estimation
methods to be compared in terms of the limit distributions
of their respective $N_{\eps,1}$ and $N_{\eps,2}$.
There are also natural connections to sequential testing
and sequential confidence sets. 
 
This paper provides general solutions to (ii) and (iii) and 
several related problems. The answer to question (1.2) turns out to be
$m\doteq2.241^2\sigma_0^2/\eps^2$, for example. 
In Section~2 the limit distribution of $\eps^2N_\eps$
is found in the i.i.d.~case, under natural conditions, 
also in the more laborious $p$-dimensional case, where 
a result is reached for a general distance function
$\|\hatt\theta_n-\theta_0\|$.
The limit distribution is that of the maximum of a certain 
squared mean zero Gau\ss ian process. 
Section 3 demonstrates that these `natural conditions' are fulfilled 
in important classes of cases, including smooth functions of averages and 
maximum likelihood estimators. 
Comparing estimators in terms of limit distributions 
for their $N_\eps$'s is seen to lead to the familiar 
expression of asymptotic relative efficiency,
that of a ratio of inverse variances, 
in the one-parameter case.   
Our arguments establish still another asymptotic optimality
property for the parametric maximum likelihood estimator sequence,
also in the $p$-parameter case:  
No other sequence will have its tail stochastically faster
included in a given $\eps$-neighbourhood, regardless of
distance measure used. 

In Section 4 a somewhat grander problem is solved,
that of obtaining a limit distribution theorem 
for the last $n$ for which the supremum distance
$\|F_n-F\|$ exceeds $\eps$, 
where $F_n$ is the empirical distribution function,
i.e., for the last $n$ in the Glivenko--Cantelli theorem.
A certain optimality property for $F_n$ is established. 
Section 5 considers the $N_\eps$ problem in a
different context, that of nonparametric density estimation. 
In this situation $\eps^2N_\eps$ goes to infinity;
it is $\eps^{5/2}N_\eps$ which has a limiting distribution.
Section 6 goes back to the situation of Sections 2 and 3, 
is technical, and demonstrates the convergence of $\eps^2\,EN_\eps$
to the appropriate limit, again under natural conditions.

In Section 7 the general methods of 
earlier sections are used to establish convergence results
for other natural quantities related to the full estimator
sequence, like $Q_\eps$, the number of times the estimator 
misses with more than $\eps$. 
Once more there is an asymptotic optimality property
for the parametric maximum likelihood method: no other estimator 
sequence has stochastically fewer $\eps$-misses. 
And a result obtained for nonparametric density estimation  
is that the best smoothing parameter for the kernel method,
in the sense of leading to the fewest $\eps$-misses 
as well as to the smallest last $n$,  
is equal to 1.008 times the traditional suggestion. 
Finally Section 8 contains a number of additional results and remarks. 


These problems have only rarely been discussed in the literature.
Bahadur (1968) considered a variable similar to our $N_\eps$, 
and indeed asked (on p.~307)
`What else can be said about $N_\eps$ [than that it is a.s.~finite],
especially for very small $\eps$?' He derived only a log-log law
for $N_\eps$, however, and failed to find 
what he [also] was searching for (p.~308): 
a simple $N_\eps$-related criterion for
comparison of estimators that would be equivalent to the 
traditional measure of asymptotic relative efficiency. 
We find such a relation, however, as mentioned above; see (2.5) and (7.2). 
Robbins, Siegmund, and Wendel (1968) found in effect 
the limit distribution of $\eps^2N_\eps$ 
for a one-sided version of the problem, 
but only in the case of a simple average of zero mean unit 
variance variables. They phrased their result in the
more probabilistic guise of the last exit time for sample 
sums $S_n$ outside the linear $n\eps$ boundary.  
Kao (1978) generalised some of their
results, and proved convergence of moments under minimal conditions,
but still only in the probabilistic random walk case just mentioned, 
which in our statistical reformulation corresponds to estimators
of the simple i.i.d.-average form. 
Some results of M\"uller (1968, 1970, 1972) 
also turn out to be related to some of ours,
as explained in Remark (i) in Section 2.  
Finally Stute (1983) proved a more statistically inspired 
result, similar to ours of Section 2, but only for certain
$M$-estimators of a simple one-dimensional location parameter. 
We emphasize that our results cover general classes of 
multi-dimensional estimators (and even an infinite-dimensional
case, in Section 4) as well as general distance measures. 

Our results give natural asymptotic relative efficiency measures 
for comparison of estimators, like (7.2) in the $p$-dimensional 
case, but are of `first order' and cannot distinguish between
competing sequences with the same first order limit distribution.  
Some `second order' results 
are in Hjort and Fenstad (1991) and Hjort (1991). 
These lead to measures of asymptotic relative deficiency 
in cases where the asymptotic relative efficiency is 1, 
and make it possible to exhibit estimator sequences that
in the second order sense 
have the smallest possible expected number of $\eps$-errors. 
See also 8E. 

\bigskip
{\bf 2. Limit distribution of $\eps^2N_\eps$.}
A simple but fundamental lemma is the following:

\smallskip
{\csc Lemma.} {{\sl
Consider i.i.d.~variables $Z_i$ with mean zero and 
variance 1, and let $S_n=\sum_{i=1}^nZ_i$. Then 
$$\sqrt{m}\sup_{n\ge m}|S_n/n|
        \rightarrow_d 
	W_{\max}=\sup_{t\ge1}|W(t)/t|
	        =_d \max_{0\le s\le1}|W(s)|, \eqno(2.1)$$
where $W(.)$ is the Brownian motion process. }} 

\smallskip
{\csc Proof:}
By Donsker's theorem, 
$S_{[mt]}/\sqrt{m}$ converges in distribution to Brownian motion $W(t)$,
a Gau\ss ian mean zero process with independent increments and
covariance function $\min(s,t)$, 
in each of the function spaces $D[b,c]$, 
see for example Billingsley (1968). Hence
$\sqrt{m}S_{[mt]}/[mt]$ tends to $W(t)/t$ in $D[1,c]$, 
from which it follows, by continuity of the supremum mapping, that
$$\sqrt{m}\sup_{m\le n\le cm}\Big|{S_n\over n}\Big|
	  =\sup_{1\le t\le c}\sqrt{m}\Big|{S_{[mt]}\over [mt]}\Big| 
	  \rightarrow_d \sup_{1\le t\le c}\Big|{W(t)\over t}\Big| 
		  =_d\sup_{c^{-1}\le s\le 1}\big|W(s)\big| $$
(employing the trick that $W^*(s)=sW(1/s)$ is a new Brownian motion).
The stronger statement of the lemma follows 
from this provided we can demonstrate 
$$\gamma_c=\limsup_{m\rightarrow\infty}{\rm Pr}\{\sqrt{m}\sup_{n\ge cm}
	 |S_n/n|\ge\delta\}\rightarrow0$$
as $c$ grows to infinity, for each given positive $\delta$; 
cf.~for example Billingsley's (1968) Theorem 4.2. 
But $\gamma_c\le6.75/c\delta^2$, as a consequence of 
a special case of inequality (6.4) 
stated and proved in Section 6. \firkant 

\smallskip
This proves useful. Assume that a one-dimensional 
$\hatt\theta_n$ admits a representation of the type
$$\hatt\theta_n-\theta_0
	 =\sigma_0\bar Z_n+R_n, \eqno(2.2)$$
where $\bar Z_n=S_n/n$ is the average of $Z_i$'s that are 
i.i.d.~with mean zero and variance 1, 
$\sigma_0$ is the standard deviation of the limiting distribution, 
and $R_n$ is the residual noise, typically of size $O_p(1/n)$. 
Define $N_\eps$ as in (1.3), 
let $y>0$ be given, and let $m$ and $y_0\doteq y$ be as in (1.4).
Then, when $\eps\rightarrow0$, 
which is the same as $m\rightarrow\infty$, 
$$\eqalign{
{\rm Pr}\{\eps^2N_\eps\ge y\}&={\rm Pr}\{N_\eps\ge m\} \cr
	&={\rm Pr}\{\sqrt{m}\sup_{n\ge m}
			 |\hatt\theta_n-\theta_0|\ge\sqrt{y_0}\} \cr
	&={\rm Pr}\{\sigma_0\sqrt{m}\sup_{n\ge m}|S_n/n
			 +\sigma_0^{-1}R_n|\ge\sqrt{y_0}\} 
 	 \rightarrow {\rm Pr}\{\sigma_0W_{\max}\ge\sqrt{y}\}, \cr}$$ 
provided the $R_n$'s are small enough. What is required is 
precisely that the difference between 
$\sqrt{m}\sup_{n\ge m}|\sigma_0S_n/n+R_n|$
and $\sigma_0\sqrt{m}\sup_{n\ge m}|S_n/n|$ goes to zero in probability 
as $m$ tends to infinity. 
For this it suffices that 
$$D_m=\sqrt{m}\sup_{n\ge m}|R_n|\rightarrow_p0, \eqno(2.3)$$
since the absolute value of the difference is dominated by $D_m$. 
(Note that the requirement for convergence of 
$\sqrt{n}(\hatt\theta_n-\theta_0)$ to $N(0,\sigma_0^2)$ 
is the weaker $\sqrt{n}R_n\rightarrow_p0$.) 
Accordingly, we have the basic result
$$\eps^2N_\eps\rightarrow_d \sigma_0^2 W_{\max}^2 \eqno(2.4)$$
for any estimator that admits representation (2.2) under condition (2.3). 
The next section demonstrates that (2.2) with (2.3) hold for
smooth functions of averages and for maximum likelihood type estimators.

\smallskip
{\csc Remarks.} 
(i) 
The Lemma was proved using just familiar 
$D[1,c]$-convergence and an extra inequality for $\gamma_c$ 
to take care of $[c,\infty)$. An alternative and in some sense 
more elegant approach is to demonstrate 
convergence in distribution of $S_{[mt]}/\sqrt{m}$ 
to $W(t)$ on the full {half}{line} $[0,\infty)$,
in some appropriate metric space of functions,
and then apply the $x(.)\rightarrow\sup_{t\ge1}|x(t)/t|$ mapping. 
One `appropriate space' is that of all 
right-continuous functions $x(.)$ with
left-hand limits satisfying $x(0)=0$ and 
$\lim_{t\rightarrow\infty}x(t)/t=0$, equipped with 
the topology induced by the norm $\sup_{t\ge0}|x(t)|/\max\{1,t\}$.
Convergence can indeed be proved using the tail inequality for $\gamma_c$,
and is also related to what is proved in M\"uller (1968);  
see also M\"uller (1970, 1972). 

(ii) 
Note that the limiting distribution in (2.4) 
is only dependent upon $\sigma_0$, and that 
the competition criterion of achieving the stochastically 
smallest limit distribution for $N_\eps$ becomes equivalent to that
of achieving the smallest possible limiting variance. 


  
(iii) 
{\sl Another optimality property for 
the maximum likelihood estimator:} 
Consider estimation in some given parametric model. 
From the previous remark it is clear that 
under traditional regularity conditions (see Section 3), 
no other sequence of estimators will have its tail 
$\{\tilda\theta_n\colon n\ge m\}$ 
included in a given neighbourhood stochastically faster than 
the sequence of maximum likelihood solutions. 
Of course the same is true for the rather wide class of estimator
sequences that are asymptotically equivalent to these, 
like Bayes estimators. 


(iv) {\sl Asymptotic relative efficiency:} 
Let $\hatt\theta_{n,1}$ and $\hatt\theta_{n,2}$ be two 
estimator sequences, both strongly consistent for $\theta_0$,
and suppose $\sqrt{n}(\hatt\theta_{n,j}-\theta_0)$ 
tends to $N(0,\sigma_j^2)$ for $j=1,2$.
The traditional notion of asymptotic relative efficiency 
is the limiting ratio of sample sizes needed by method 1 and method 2 
to achieve some desirable accuracy.  
The classical formula for the a.r.e.~of method 2 w.r.t.~method 1 is 
${\rm a.r.e.}=\sigma_1^2/\sigma_2^2$, 
motivated from approximate risk or from approximate
length of confidence intervals, see Serfling (1980, p.~50--52). 
This measure also plays a natural r\^ole in the Pitman
approach to comparing test statistics, see Serfling's Section 10.2.    
Define last $n$ variables $N_{\eps,1}$ and $N_{\eps,2}$ 
for the two estimator sequences. 
The ratio of sample sizes viewpoint invites the following 
as natural measures of asymptotic relative efficiency:    
$$\lim_{\eps\rightarrow0}
	{{\rm med}\{N_{\eps,1}\}\over {\rm med}\{N_{\eps,2}\}}
	={\sigma_1^2\over \sigma_2^2}={\rm a.r.e.}, \quad 
  \lim_{\eps\rightarrow0}
	{EN_{\eps,1}\over EN_{\eps,2}}
	={\sigma_1^2\over \sigma_2^2}={\rm a.r.e.} \eqno(2.5)$$  
(see Section 6 for convergence of moments). 
This provides fresh and independent motivation for the a.r.e.~measure. 
See also 7A, 7B, and 8E. 

\smallskip
Let us now turn to the $p$-dimensional case.
Let $N_\eps$ be defined as in (1.3) but with respect to some given
distance function $\|\hatt\theta_n-\theta_0\|$ in ${\cal R}^p$,
for example ordinary Euclidian distance. 
We have primarily distances of the type $\{(x-y)'A(x-y)\}^{1/2}$ in mind,
where $A$ is symmetric and positive definite, 
but require only that
$\|x\|$ is a function on $p$-vectors with the properties 
$\|x+y\|\le\|x\|+\|y\|$, 
$\|x\|=0$ if and only if $x=0$, 
$\|x_n\|\rightarrow\|x\|$ when $x_n\rightarrow x$,
$\|ax\|=|a|\,\|x\|$ for scalars $a$,  
and $\|x\|=\|(x_1,\ldots,x_p)'\|\le c\sum_{i=1}^p|x_i|$ for
some constant $c$. See 8D for other distances. 

\smallskip
{\csc Theorem.}  
{{\sl Suppose that 
$$\hatt\theta_n-\theta_0=\Sigma_0^{1/2}
	 {1\over n}\sum_{i=1}^nZ_i+R_n, \eqno(2.6)$$
where the $Z_i$'s are i.i.d.~with zero mean and the $p\times p$
identity matrix as covariance matrix. Suppose further that 
$$D_m=\sqrt{m}\sup_{n\ge m}\|R_n\|\rightarrow_p0; \eqno(2.7)$$
in particular $\sqrt{n}(\hatt\theta_n-\theta_0)\rightarrow_dN_p(0,\Sigma_0)$. 
Let $G_p(s)=\Sigma_0^{1/2}W(s)$, where 
$W(s)=(W_1(s),\allowbreak\ldots,W_p(s))'$ is a vector of $p$ independent 
Brownian motions, each evaluated at the same $s$. 
Then, as $\eps$ tends to zero, }}
$$\eps^2N_\eps\rightarrow_d G_{p,\max}^2
        =\bigl\{\sup_{0\le s\le1}\|G_p(s)\|\bigr\}^2
	 =\sup_{0\le s\le1}\|\Sigma_0^{1/2}W(s)\|^2. \eqno(2.8)$$

\smallskip
{\csc Proof:} Somewhat more elaborate arguments are necessary now.
We prove first that 
$$\sqrt{m}\sup_{n\ge m}\|\Sigma_0^{1/2}S_n/n\|
	 \rightarrow_d G_{p,\max}
	 =\sup_{0\le s\le1}\|\Sigma_0^{1/2}W(s)\|, \eqno(2.9)$$
where again the $S_n$'s are partial sums of the $Z_i$'s. 
Observe first that the stochastic process 
$(S_{1,[mt_1]}/\sqrt{m},\ldots,S_{p,[mt_p]}/\sqrt{m})'$,
where $S_{j,n}$ is the $j$'th component of $S_n$, converges 
to $(W_1(t_1),\ldots,W_p(t_p))'$ in each $D[b,c]^p$,
equipped with the product Skorohod topology. 
[This $p$-variate version of Donsker's theorem follows from 
the 1-variate theorem by tightness and finite-dimensional convergence.]
By the continuous mapping theorem, 
$\Sigma_0^{1/2}\sqrt{m}S_{[mt]}/[mt]$ converges to $\Sigma_0^{1/2}W(t)/t$
in $D_p[1,c]$, the space of all right-continuous functions
$[1,c]\rightarrow{\cal R}^p$ with left hand limits, equipped with the 
Skorohod topology. And since the supremum mapping is continuous too, 
$$\eqalign{
\sqrt{m}\sup_{m\le n\le cm}\|\Sigma_0^{1/2}S_n/n\|
	 &=\sqrt{m}\sup_{1\le t\le c}\|\Sigma_0^{1/2}S_{[mt]}/[mt]\| \cr
	 &\rightarrow_d\sup_{1\le t\le c}\|\Sigma_0^{1/2}W(t)/t\| 
	  =_d\sup_{c^{-1}\le s\le1}\|\Sigma_0^{1/2}W(s)\|. \cr}$$
Claim (2.9) follows since $\gamma_c=\limsup_{m\rightarrow\infty}
{\rm Pr}\{\sqrt{m}\sup_{n\ge cm}\|\Sigma_0^{1/2}S_n/n\|\ge\delta\}$
tends to zero as $c$ grows, by a simple inequality relating this quantity 
to a sum of $p$ one-dimensional analogues; cf.~the proof of the Lemma.

The rest of the proof follows from (2.9) and regularity condition (2.7).
For let again $m$ and $y_0$ be as in (1.4). Then 
$${\rm Pr}\{\eps^2N_\eps\ge y\}
	={\rm Pr}\{\sqrt{m}\sup_{n\ge m}\|\Sigma_0^{1/2}S_n/n
			 +R_n\|\ge\sqrt{y_0}\} $$
is seen to converge to 
${\rm Pr}\{\sup_{0\le s\le 1}\|\Sigma_0^{1/2}W(s)\|\ge\sqrt{y}\}$, 
which is the same as 
$\eps N_{\eps}^{1/2}\rightarrow_d G_{p,\max}$ 
or $\eps^2N_\eps\rightarrow_d G_{p,\max}^2$. \firkant 

\smallskip
In a parametric model 
the maximum likelihood estimator sequence achieves the
smallest possible limit covariance matrix and therefore also achieves 
the stochastically smallest possible limit distribution for $N_\eps$,
regardless of distance measure $\|\hatt\theta_n-\theta_0\|$, 
cf.~Remark (iii) above for the one-dimensional case.

\smallskip
{\csc Corollary.} {{\sl 
Let conditions be as in the theorem, and let  
$\|\hatt\theta_n-\theta_0\|=
\{(\hatt\theta_n-\theta_0)'\Sigma_0^{-1}(\hatt\theta_n-\theta_0)\}^{1/2}$ 
be $\Sigma_0$-weighted Mahalanobis distance. Then 
$$\sqrt{m}\sup_{n\ge m}\|\hatt\theta_n-\theta_0\|
	\rightarrow_d\chi_{p,\max} 
	\quad {\sl and} \quad
	\eps^2N_\eps\rightarrow_d \chi_{p,\max}^2, \eqno(2.10)$$
as respectively $m\rightarrow\infty$ and $\eps\rightarrow0$, 
where $\chi_{p,\max}^2=\max_{0\le s\le1}\sum_{i=1}^pW_i(s)^2$. }} 

\smallskip
In particular the limit distribution is the very same one 
in each estimation problem with $p$ parameters! 
One can prove that (2.10) continues to hold
when $\Sigma_0$ is replaced by a strongly consistent estimate $\hatt\Sigma_n$ 
($\hatt\Sigma_n\rightarrow_p\Sigma_0$ does not suffice).
Details are in Hjort and Fenstad (1990). 

The in some sense natural extension of (2.5) to the $p$-dimensional
case would be ${\rm a.r.e.}=EH(\Sigma_1)/EH(\Sigma_2)$, where 
$H(\Sigma)=\max_{0\le s\le 1}W(s)'\Sigma W(s)$,
since $\eps^2EN_\eps\rightarrow EH(\Sigma)$ under Euclidian distance. 
There is no simple formula for $EH(\Sigma)$, however. 
A simple explicit a.r.e.~measure emerges in 7B. 


\bigskip
{\bf 3. Applications to special cases.}
In this section we confirm that the necessary regularity conditions
(2.6) and (2.7) indeed pertain in the usual situations,
both under parametric and nonparametric circumstances. 
In the one-dimensional case, suppose 
that $\hatt\theta_n$ admits the representation
$$\hatt\theta_n-\theta_0
	 ={1\over n}\sum_{i=1}^n\sigma_0Z_i
	  +R_n, \quad R_n=\delta_n\,\bar U_n
	 =\delta_n{1\over n}\sum_{i=1}^nU_i, \eqno(3.1)$$ 
where the $U_i$'s are i.i.d.~with mean zero and finite variance, and 
$\delta_n\rightarrow0$ a.s. Then (2.3) holds, since
$$D_m=\sqrt{m}\sup_{n\ge m}|\delta_n\bar U_n|
	 \le \sup_{n\ge m}|\delta_n|\,
	 \sqrt{m}\sup_{n\ge m}|\bar U_n|\rightarrow_p0,$$
in that the second term has a limit in distribution, 
by the Lemma of Section 1,
and the first term tends to zero in probability, by the definition of 
$\delta_n\rightarrow0$ a.s. There is a similar result for the 
$p$-dimensional case: 
If (3.1) holds, with $\Sigma_0$ replacing $\sigma_0^2$, 
and where the $U_i$'s are i.i.d.~vectors
with mean zero and finite covariance matrix, and $\delta_n$ is a matrix
with components that all tend to zero a.s., then $D_m$ of (2.7)
tends to zero in probability. 
This follows essentially by the one-dimensional
argument. To see this, let ${\rm norm}(\delta_n)$
be the matrix norm of $\delta_n$, defined as the maximum 
of $\|\delta_n x\|$ over $\|x\|\le 1$;
for Euclidian distance-norm ${\rm norm}(\delta_n)$
is equal to the largest eigen-value, for example. 
Then $\|\delta_n\bar U_n\|\le {\rm norm}(\delta_n)\|\bar U_n\|$,
which goes a.s.~to zero by the continuity of the $\|\cdot\|$ norm.  

\smallskip 
{\sl 3A. Smooth functions of averages.} 
Suppose $\hatt\theta_n=h(\bar B_n)$
and $\theta_0=h(b)$, where $\bar B_n$ is the average of i.i.d.~variables
$B_i$ with $EB_i=b$ and ${\rm Var}\,B_i=\tau^2$. 
If $h$ has a continuous derivative in a neigbourhood of $b$, then  
$$\sqrt{n}(\hatt\theta_n-\theta_0)=h'(b)\sqrt{n}(\bar B_n-b)
	  +\{h'(\tilda b_n)-h'(b)\}\sqrt{n}(\bar B_n-b)$$
for some random $\tilda b_n$ between $b$ and $\bar B_n$. 
This is as in (3.1) with $\delta_n=h'(\tilda b_n)-h'(b)$. 
But it is easy to see that $\delta_n\rightarrow0$ 
a.s.~by the strong law of large numbers for $\bar B_n$. 
Hence (2.4) holds, with $\sigma_0^2=h'(b)^2\tau^2$. --- More generally, 
suppose $\hatt\theta_n$ is $p$-dimensional and that 
$\hatt\theta_j=h_j(\bar B_{n,1},\ldots,\bar B_{n,r})$
for $j=1,\ldots,p$,  
for $r$ averages of i.i.d.~vectors $(B_{i,1},\ldots,B_{i,r})'$
with mean $b=(b_1,\ldots,b_r)'$ and finite covariance matrix $T$,
and let $h_j(b)=\theta_{0,j}$.  
If only $h_1,\ldots,h_p$ have Jacobi matrix $J(x)$ 
with partial derivatives $\dell h_j(x)/\dell x_l$ 
that are continuous in a neighbourhood of $(b_1,\ldots,b_r)$,
then (2.7) holds. And this implies (2.8) with 
$\Sigma_0=J(b)TJ(b)'$. 

\smallskip
{\csc Example 1.} 
Suppose $X_1,X_2,\ldots$ are i.i.d.~with finite sixth moment.
Then $\hatt\theta_n={1\over n}\sum_{i=1}^n(X_i-\bar X)^3$,
the natural and strongly consistent estimator of $\theta_0=E(X_i-EX_i)^3$,
is a smooth function of the sample averages of $X_i$, $X_i^2$, $X_i^3$,
and $\eps^2N_\eps\rightarrow_d\sigma_0^2W_{\max}^2$, 
where $\sigma_0^2$ is the limit variance of 
$\sqrt{n}(\hatt\theta_n-\theta_0)$.
[In fact $\sigma_0^2=(9+\alpha_6-6\alpha_4-\alpha_3^2)\tau^6$,
where $\alpha_p=E(X_i-EX_i)^p/\tau^p$ and $\tau$ is the standard deviation
for $X_i$.]

\smallskip
{\sl 3B. Maximum likelihood estimators.} 
The typical argument that leads to a limit distribution result for
the maximum likelihood estimator uses Taylor expansion to get 
$$\sqrt{n}(\hatt\theta_n-\theta_0)
	  =J_n^{-1}{1\over \sqrt{n}}\sum_{i=1}^nU(X_i)
	  \rightarrow_d N_p\{0,J_0^{-1}\}, \eqno(3.2)$$
where $U(X_i)=\dell\log f(X_i,\theta_0)/\dell\theta$ 
is the score function and $J_n\rightarrow_pJ_0$, 
the variance matrix for $U(X_i)$
computed under $f(x,\theta_0)$, 
i.e.~the familiar Fisher information matrix. 
It follows from previous arguments that (3.2) also secures 
convergence in distribution of $\eps^2N_\eps$, 
as in (2.8), with $\Sigma_0=J_0^{-1}$, 
provided there also is a.s.~convergence $J_n\rightarrow J_0$. 
But this is true, under weak conditions. 
It is for example not difficult to
prove that the conditions used in Lehmann's (1983) Section 6.4
suffice. 

One can also prove that if the model specifies $f(x,\theta)$,
but the true density $f$ does not belong, then (2.8) holds again
under mild regularity conditions, 
but with a different interpretation of $\theta_0$ and a different matrix. 
The $\theta_0$ that now enters is not `true', but rather `least false'
or `best fitting', and can be characterised as the parameter value
that minimises the Kullback--Leibler distance 
$d[f\colon f(.,\theta)]=\int f(x)\log\{f(x)/f(x,\theta)\}\,dx$. 
Furthermore, $\Sigma_0=J_0^{-1}K_0J_0^{-1}$, 
where $K_0$ is the variance matrix, under the true $f$, 
of the score function computed at $\theta_0$; 
and $J_0$ is minus the expected value, under the true $f$, 
of the twice differentiated log-density also computed at $\theta_0$. 
If the model happens to be perfect, 
then $\theta_0$ deserves to be called `true', and $J_0=K_0$. 
Proofs and discussion of these claims about maximum likelihood 
under the agnostic viewpoint can be found in Hjort (1986, Ch.~3). 


\smallskip
{\csc Example 2.} 
Consider again maximum likelihood estimation in 
a given parametric family $f(x,\theta)$. 
Let distance be measured in the invariant Mahalanobis way, 
$\|\theta-\theta_0\|^2=(\theta-\theta_0)'J_0(\theta-\theta_0)$,
and let $N_\eps$ be the last $n$ for which 
$\|\hatt\theta_n-\theta_0\|\ge\eps$. 
Then $\eps^2N_\eps$ tends to 
$\max_{0\le s\le 1}W(s)'J_0^{-1/2}K_0J_0^{-1/2}W(s)$
(which is $\chi^2_{p,\max}$ if the model is correct).
For a specific example, suppose the model specifies the normal
density $f(x,\theta)=f(x,\mu,\sigma)$, but assume only that 
the true $f$ is symmetric with finite fourth moment.
Then the least false parameters are 
$\mu_0=E_fX_i$ and $\sigma_0={\rm stdev}_fX_i$. 
One also finds 
$J_0^{-1/2}K_0J_0^{-1/2}={\rm diag}(1,1+\half\beta_2)$,
where $\beta_2=E\{(X-\mu_0)/\sigma_0\}^4-3$ is the kurtosis.
Hence $\eps^2N_\eps$ tends to 
$\max_{0\le s\le 1}\{W_1(s)^2+(1+\half\beta_2)W_2(s)^2\}$,
where $W_1$ and $W_2$ are independent Brownian motions.

\smallskip
{\csc Example 3.}
Let $(Y_1,\ldots,Y_p)$ be multinomial $(n,\theta_1,\ldots,\theta_p)$,
with $\sum_{i=1}^p\theta_i=1$ and $\sum_{i=1}^pY_i=n$. 
Let $N_\eps$ be the last $n$ at which 
$\sum_{i=1}^p(\hatt\theta_i-\theta_i)^2/\theta_i\ge\eps^2$,
where $\hatt\theta_i=Y_i/n$ is the usual maximum likelihood
estimator of $\theta_i$. 
This corresponds in fact to measuring distance from 
$(\hatt\theta_1,\ldots,\hatt\theta_{p-1})$ to 
$(\theta_1,\ldots,\theta_{p-1})$ in the Mahalanobis way,
see the Corollary ending Section 2. 
Hence $\eps^2N_\eps$ tends to $\chi^2_{p-1,\max}$. 
The same is true if $\theta_i$'s are replaced with $\hatt\theta_i$'s
in the denominators, see Section 8F. 


\smallskip
{\sl 3C. Differentiable functionals.} 
In many situations the estimator $\hatt\theta_n$ 
can be thought of as a functional $T$ evaluated at 
the empirical distribution function $F_n$,
while the true parameter $\theta_0$ correspondingly 
is equal to $T(F)$ for the true $F$. 
Suppose $T$ is so-called locally Lipschitz differentiable at $F$
w.r.t.~the supremum norm $\|G-F\|=\sup_x|G(x)-F(x)|$,
which means that 
$T(G)-T(F)=\int I(F,x)\{dG(x)-dF(x)\}+O(\|G-F\|^2)$, 
featuring the influence function 
$I(F,x)=\lim_{\eps\rightarrow0}\{T((1-\eps)F+\eps\delta_x)-T(F)\}/\eps$.
This might be interpreted as a reasonable minimum amount of smoothness
on the part of $T(.)$. Examples are given in Shao (1989),
including general $L$- and $M$-estimators. In particular
the somewhat non-smooth median functional is still 
locally Lipschitz differentiable. --- Under this assumption 
it holds that 
$$\hatt\theta_n-\theta_0={1\over n}\sum_{i=1}^nI(F,X_i)
	 +O(\|F_n-F\|^2).$$
But it is known that 
$\|F_n-F\|^2\le Kn^{-1}\log\log n$ a.s., for some large $K$,
see for example Shao (1989). It follows that 
$|R_n|\le K'n^{-1}\log\log n$ in representation (2.2), 
and this implies (2.3). Consequently (2.4), 
or (2.8) in the $p$-dimensional case, 
are true for functionals that are locally Lipschitz differentiable.  

\bigskip
{\bf 4. The last $n$ for Glivenko--Cantelli.}
Let $X_1,X_2,\ldots$ be independent from some continuous $F$, 
and let $F_n(t)$ 
be the empirical distribution function 
${1\over n}\sum_{i=1}^nI\{X_i\le t\}$ 
based on the first $n$ data points. Then 
$$\eqalign{
\sqrt{m}\sup_{n\ge m}|F_n(t)-F(t)|
	&\rightarrow_d\{F(t)(1-F(t))\}^{1/2}\,W_{\max}, \cr
	\eps^2N_\eps(t)&\rightarrow_dF(t)(1-F(t))\,W_{\max}^2 \cr} \eqno(4.1)$$
by previous efforts,  
where $N_\eps(t)$ is the last $n$ for which $|F_n(t)-F(t)|\ge \eps$.
Can we obtain similar results for the supremum 
distance $\|F_n-F\|$? 

The answer to these somewhat grander questions 
must involve asymptotic arguments in $n$ and $t$
simultaneously. Let $K_0(s,t)$ be a Kiefer process on 
$[0,\infty)\times[0,1]$. This is a two-parameter 
zero mean Gau\ss ian process with continuous sheets and 
$${\rm cov}\{K_0(s_1,t_1),K_0(s_2,t_2)\}
	=(s_1\wedge s_2)\,(t_1\wedge t_2-t_1t_2). \eqno(4.2)$$
It behaves like a Brownian bridge in $t$ for fixed $s$ and 
like Brownian motion in $s$ for fixed $t$. Note that 
$K(s,t)=sK_0(s^{-1},t)$ is another Kiefer. 

\smallskip
{\csc Theorem.} {{\sl
Let $N_\eps$ be the last $n$ at which $\|F_n-F\|\ge\eps$, and 
let $K_{\max}$ be the maximum of $|K(s,t)|$ over the 
unit square $[0,1]\times[0,1]$. Then 
$$\sqrt{m}\sup_{n\ge m}\|F_n-F\|\rightarrow_d K_{\max}
	\quad {\sl and} \quad
	\eps^2N_\eps\rightarrow_d K_{\max}^2,$$
as respectively $m\rightarrow\infty$ and $\eps\rightarrow0$.}}

\eject 

\smallskip
{\csc Proof:} 
Considerations involving the inverse transformation 
$X_i'=F^{-1}(\xi_i)$, where the $\xi_i$'s are i.i.d.~from 
the uniform distribution $F_0(t)=t$ 
on $[0,1]$, reveal that the distribution of the full sequence
of $\|F_n-F\|$ is equal to that of $\|F_{n,0}-F_0\|$,
where $F_{n,0}$ is the empirical distribution of the $n$
first $\xi_i$'s. Accordingly we might as well take  
$F$ to be $F_0$ from the outset, and this simplifies matters below. 

The LeCam--Bickel--Wichura--M\"uller theorem states that 
the process 
$$K_m(s,t)={1\over \sqrt{m}}\sum_{i=1}^{[ms]}
	\bigl[I\{X_i\le t\}-F_0(t)\bigr]
	={[ms]\over \sqrt{m}}\bigl\{F_{[ms]}(t)-t \bigr\}$$
converges in distribution to $K_0(s,t)$ in $D\{[b,c]\times[0,1]\}$
with the Skorohod metric, for each $[b,c]$ interval, 
see for example Shorack and Wellner (1986, Chapter 3.5). 
For us it is more convenient to study 
$$H_m(s,t)=\sqrt{m}\bigl\{F_{[ms]}(t)-t \bigr\}
	={m\over [ms]}K_m(s,t)
	\rightarrow_d (1/s)K_0(s,t)=K(1/s,t).$$
(The (4.1) results follow anew from this.) 
By the continuous mapping theorem 
$$\sqrt{m}\sup_{n\ge m}\|F_n-F_0\|
	=\sup_{s\ge 1}\sup_{0\le t\le 1} |H_m(s,t)| 
	 \rightarrow_d\sup_{s\ge 1}\sup_{0\le t\le 1}|K(1/s,t)|
	 	=K_{\max}. $$
Reasoning once more as in (1.4) we also obtain 
$\eps^2N_\eps\rightarrow_dK_{\max}^2$. 
This incidentally also gives a sequential fixed-width 
nonparametric simultaneous confidence band for $F$. 

The argument presented here is heuristic at one point, 
since convergence in distribution of the $H_m$ process is only 
guaranteed on each $[1,c]\times[0,1]$. Therefore
only convergence of 
$\sqrt{m}\sup_{m\le n\le cm}\|F_n-F\|$ to the maximum of 
$|K(s,t)|$ over $[1/c,1]\times[0,1]$ is rigorously proved, so far.
What needs to be ascertained is that 
$$\gamma_c=\limsup_{m\rightarrow\infty}{\rm Pr}\{\sqrt{m}
	\sup_{n\ge cm}\|F_n-F_0\|\ge \delta\}\rightarrow0
	\quad{\rm as\ }c\rightarrow\infty, \eqno(4.3)$$
cf.~once more Billingsley's (1968) Theorem 4.2
and the corresponding technical point in the proof of Lemma of Section 2.
It will suffice to prove 
$${\rm Pr}\{\sqrt{m}\sup_{n\ge m}\|F_n-F_0\|\ge b\}\le A/b^4
	\quad {\rm for\ all\ } b {\rm \ and\ }m, \eqno(4.4)$$
for some large enough constant $A$, 
since this implies $\gamma_c\le A/c^2\delta^4$.  
[An alternative route is to prove weak convergence 
in some appropriate function space on $[1,\infty)\times[0,1]$,
with suitable metric. A result of M\"uller (1970) is of this type.]
 
To prove (4.4) we shall use general fluctuation inequalities
provided by Bickel and Wichura (1971) for two-parameter
processes. For neighbouring blocks $B$ and $C$ in the unit 
square one can show
$E\{K_m(B)^2K_m(C)^2\}\le 3\mu(B)\mu(C)$,
where $\mu$ is Lebesgue measure, see Shorack and Wellner
(1986, Chapter 3.5). This implies 
${\rm Pr}\{\sup_{s,t\in[0,1]}|K_m(s,t)|\ge b\}\le A/b^4$,
for some universal constant $A$, by Bickel and Wichura's 
Theorem 1 in conjunction with their inequality (1). 
But 
$$\eqalign{
\half\sqrt{m}\max_{m/2\le n\le m}\|F_n-F_0\|
	&\le \sqrt{m}\max_{m/2\le n\le m}
		\hbox{$n\over m$}\|F_n-F_0\| \cr
	&\le \sqrt{m}\max_{n\le m}\hbox{$n\over m$}\|F_n-F_0\|
	 =\sup_{s,t\in[0,1]}|K_m(s,t)|. \cr}$$
This is soon translated into 
${\rm Pr}\{\sqrt{m}\max_{m\le n\le 2m}\|F_n-F_0\|\ge b\}
	\le A/(b/\sqrt{2})^4$,
for all $m$ and $b$. Let $2^k\le m<2^{k+1}$. Then
the left-hand side of (4.4) is bounded by the sum of 
${\rm Pr}\{\sqrt{2^i}\max_{2^i\le n<2^{i+1}}\|F_n-F_0\|\ge 
	\sqrt{2^i}b/\sqrt{m}\}$,
for $i\ge k$. Bounding each of these in the way just described 
gives at the end of the night (4.4), with constant ${64\over 3}A$,
which concludes our proof. \firkant 


\smallskip
The two-parameter stochastic process approach is very powerful,
and allows us to reach other related results as well. 
As but one example, let ${\rm CM}_n^2=\int\{F_n(t)-F(t)\}^2\,dF(t)$
be the Cram\'er--von Mises statistic. Using the $H_m$ process 
from the proof above we have
$$m\sup_{n\ge m}{\rm CM}_n^2
	=_d \sup_{s\ge1}\int_0^1H_m(s,t)^2\,dt
	\rightarrow_d\sup_{0\le s\le 1}\int_0^1K(s,t)^2\,dt=\Lambda^2 .$$
We also have $\eps^2N_\eps\rightarrow\Lambda^2$, if
$N_\eps$ is the last $n$ where ${\rm CM}_n\ge\eps$. 
[Ways of simulating the distributions of $\Lambda^2$ 
and $K_{\max}$ are described in Hjort and Fenstad (1990).]

{\sl Asymptotic optimality of $F_n$:} Are there estimators better than
$F_n$, as measured by expected smallness of $N_\eps$,
as $\eps$ tends to zero? The answer to this
question is no, if one rules out the superefficiency phenomenon.
This follows from the H\'ajek convolution type representation
theorem for limit distributions for $\sqrt{n}(\tilda F_n-F)$ 
proved by Beran (1977), in conjunction with the arguments used above. 

\bigskip
{\bf 5. The last $n$ for nonparametric density estimators.}
Consider a kernel type estimator 
$f_n(x)={1\over n}\sum_{i=1}^nK((x-X_i)/h_n)/h_n$ 
for the unknown density $f(x)$ based on
the first $n$ data points in an i.i.d.~sequence. 
What is the size of $N_\eps$, the last time $|f_n(x)-f(x)|\ge\eps$? 
Techniques from Sections 2 and 3 can be employed to reach a
limit distribution result though some extra care is needed since 
$h_n$ varies with sample size 
(the minimum requirement for strong consistency is 
$h_n\rightarrow0$ and $nh_n\rightarrow\infty$). 

Suppose $f$ has two continuous derivatives around the
given $x$, and let the kernel density $K$ have mean zero
and unit variance and finite $\beta_K=\int K(u)^2\,du$.
Let $h_n=cn^{-1/5}$, well known to be the optimal rate. 
Study $Z_m(t)=m^{2/5}\{f_{[mt]}(x)-f(x)\}$ for the fixed $x$. 
It splits into a bias term $b_m(t)$ and a zero mean term $Z_m^0(t)$.
The first can be seen to 
converge to the function $\half c^2f''(x)/t^{2/5}$, uniformly
over finite $t$-intervals. 
The second can be proved to converge
in distribution to a Gau\ss ian zero mean process
with covariance function of the form $c^{-1}g(s/t)f(x)/t^{4/5}$,
where in fact $g(z)=z^{1/5}\int K(u)K(z^{1/5}u)\,du$. 
Hence 
$$Z_m(t)\rightarrow_d Z(t)
	=\bigl[\half c^2f''(x)+c^{-1/2}f(x)^{1/2}V(t)\bigr]/t^{2/5}$$
for a certain normal zero mean stochastic process 
$V(.)$ with constant variance $\beta_K$.
From this result, using arguments parallel to those of Section 2,
it is not difficult to derive 
$$\eps^{5/2}N_\eps\rightarrow_d 
	Z_{\max}^{5/2}=\bigl\{\sup_{t\ge1}|Z(t)|\bigr\}^{5/2}.$$ 
\qquad 
How should $c$ in $cn^{-1/5}$ be chosen? 
The approximate mean squared error is 
${1\over4}h^4f''(x)^2\allowbreak+\beta_Kf(x)/nh$, and is minimised 
for $h_n=c_0(x)n^{-1/5}$, where $c_0(x)=\{\beta_Kf(x)/f''(x)^2\}^{1/5}$.
One version of the variable kernel approach to density estimation
is to aim for this value, using a smooth pilot estimate to reach
$\hatt c_0(x)$, say. We could try to make $EN_\eps$ as small
as possible by making $E|Z(t)|^{5/2}$ as small as possible. 
But this expectation can be written 
$a^{-5/4}E|\half a^{5/4}+N(0,1)|^{5/2}$ times other terms 
not depending on $c$, where $a=c/c_0(x)$. Careful numerical integrations
reveal that minimum occurs for $a_0=1.008$. 
Hence $1.008\,c_0(x)n^{-1/5}$ is best from the $EN_\eps$ 
point of view. See also 7D. 

The derivative of $f$ is even more difficult to estimate 
with good precision. This is reflected in high values for
$N'_\eps$, the last $n$ for which $|f_n'(x)-f'(x)|$ exceeds $\eps$.
By techniques similar to those
sketched above one can prove that $\eps^{7/2}N'_\eps$ 
tends to some appropriate $(Z'_{\max})^{7/2}$ in distribution.

It would be interesting to reach results for $N_\eps$'s 
connected to global deviance 
measures like $\int(f_n-f)^2/f\,dx$ or the statistically natural
but technically difficult $\int|f_n-f|\,dx$ as well. 
Techniques from Bickel and Rosenblatt (1973) would be 
appropriate, but we haven't pursued this. 
 
\bigskip
{\bf 6. Convergence of moments.}
We have proved that $\eps^2N_\eps\rightarrow_d\sigma_0^2W_{\max}^2$
(in the one-dimensional case), 
and it is clear that $\eps^2EN_\eps$ should tend to 
$\sigma_0^2EW_{\max}^2$ under conditions pertaining to uniform
integrability. The present section derives this and a couple of 
related results under natural conditions. 

We should like to prove 
$$E\eps^2N_\eps=\int_0^\infty{\rm Pr}\{\eps^2N_\eps\ge y\}\,dy
	  \rightarrow\int_0^\infty{\rm Pr}\{\sigma_0^2W_{\max}^2\ge y\}\,dy
	  =E\sigma_0^2W_{\max}^2,$$
and this holds by Lebesgue's theorem on dominated convergence 
provided we can bound 
$${\rm Pr}\{\eps^2N_\eps\ge y\}
	  ={\rm Pr}\{\sqrt{m}\sup_{n\ge m}
	 |\hatt\theta_n-\theta_0|\ge\sqrt{y_0}\},
	  \quad m=\langle y/\eps^2\rangle,
	  \quad y_0=m\eps^2, \eqno(6.1)$$ 
with some integrable function, uniformly in $\eps$.   
A sufficient condition is therefore that
for some positive $\eps_0$ 
$${\rm Pr}\{\sqrt{m}\sup_{n\ge m}|\hatt\theta_n-\theta_0|\ge a\}
	  \le K/a^{2+\lambda} \quad 
		  {\rm \ when\ }0<a/\sqrt{m}\le\eps_0, \eqno(6.2)$$
for some positive $\lambda$ and some companion constant $K$. 

We start with the simplest case $\hatt\theta_n-\theta_0=\sigma_0S_n/n$, 
with partial sums of $Z_i$'s that are i.i.d.~with mean zero and variance 1,
as in the Lemma of Section 2. 

\smallskip
{\csc Lemma.} {{\sl 
Suppose $E|Z_i|^{2+\lambda}<\infty$ for some $\lambda\ge0$. 
Then there is a constant $c_{2+\lambda}$ such that 
$$E|S_n|^{2+\lambda}\le c_{2+\lambda}n^{1+\lambda/2}E|N(0,1)|^{2+\lambda}
	  \quad{\sl for\ all\ }n \eqno(6.3)$$
(and $c_{2+\lambda}$ can be replaced with 1.001 if we change 
`for all $n$' to `for all large $n$'). 
Furthermore, 
$${\rm Pr}\{\sqrt{m}\sup_{n\ge m}|S_n/n|\ge a\}
	  \le{6.75\,c_{2+\lambda}
		 E|N(0,1)|^{2+\lambda}\over a^{2+\lambda}}
	 \quad {\sl for\ all\ }m{\sl \ and\ }a. \eqno(6.4)$$ }}

{\csc Proof:} Of course $S_n/\sqrt{n}\rightarrow_d N(0,1)$. 
Results from von Bahr (1965) can be used to show 
$E|S_n/\sqrt{n}|^{2+\lambda}=E|N(0,1)|^{2+\lambda}+r_n$, 
where $|r_n|\le M/\sqrt{n}$
for some $M$. In particular there is convergence, and (6.3)
(with accompanying parenthetical remark) follows from this. 
As a step in the rest of the proof we utilise a generalisation of 
Kolmogorov's inequality, namely 
$${\rm Pr}\{\max_{i\le n}|S_i|\ge a\}\le E|S_n|^{2+\lambda}/a^{2+\lambda},$$
which can be found in e.g.~Lo\`eve (1960, p.~263). 
Let $q>1$, suppose $q^k\le m<q^{k+1}$, and let us abbreviate
$c_{2+\lambda}E|N(0,1)|^{2+\lambda}$ with $K$. Then 
$$\eqalign{
{\rm Pr}\{\sqrt{m}\sup_{n\ge m}|S_n/n|\ge a\}
	 &\le\sum_{i=k}^\infty{\rm Pr}\{\max_{q^i\le n<q^{i+1}}|S_n|
		 \ge aq^i/\sqrt{m}\} \cr
	 &\le\sum_{i=k}^\infty{K(q^{i+1})^{1+\lambda/2}
		 \over (aq^i/\sqrt{m})^{2+\lambda}} \cr
	 &={K\over a^{2+\lambda}}m^{1+\lambda/2}q^{1+\lambda/2}
		 \sum_{i=k}^\infty\Bigl({1\over q}\Bigr)^{i(1+\lambda/2)} 
	  \le{K\over a^{2+\lambda}}
		 {q^{3+3\lambda/2}\over q^{1+\lambda/2}-1}. \cr}$$
The best value of $q$ corresponds to 
$q^{1+\lambda/2}={3\over2}$, and the result follows. \firkant


\smallskip
Note that the right hand side of (6.4) becomes $6.75/a^2$ for
$\lambda=0$; this was needed in the proof of Section 2's Lemma.
Robbins, Siegmund, and Wendel (1968) have inequality (6.4) 
for this simplest $\lambda=0$ case 
(but with constant 8 instead of 6.75). 

This basic lemma can now be used to prove 
$E\eps^2N_\eps\rightarrow1.832\,\sigma_0^2$
in various situations. Consider smooth functions of averages.
Suppose $\hatt\theta_n=h(\bar B_n)$ estimates $\theta_0=h(b)$,
where $\bar B_n$ is the average of i.i.d.~variables $B_i$ with
mean $b$ and variance $\tau^2$, as in 3A. 
In particular $\eps^2N_\eps\rightarrow_d\sigma_0^2W_{\max}^2$,
where $\sigma_0^2=h'(b)^2\tau^2$, 
if only $h$ has a continuous derivative around $b$. 

\smallskip {{\sl
{\csc Theorem.} Suppose in addition that $E|B_i|^{2+\lambda}$ 
is finite for some positive $\lambda$. 
Then $\eps^2EN_\eps\rightarrow2G\sigma_0^2$,
where $G=0.915966 ...$ is the Catalanian constant (see 8A). }}

\smallskip
{\csc Proof:} 
We are to prove (6.2). 
This is very immediate if $h'$ is bounded, but some care is needed
to cover all the interesting cases where $h'$ is unbounded,
cf.~Example 1 of Section 3A. 
With notation as in 3A we have
$$\eqalign{
{\rm Pr}\{\sqrt{m}\sup_{n\ge m}|\hatt\theta_n-\theta_0|\ge 2a\}
	 &\le{\rm Pr}\{\sqrt{m}\sup_{n\ge m}|h'(b)(\bar B_n-b)|\ge a\} \cr
	 &\quad +{\rm Pr}\{\sqrt{m}\sup_{n\ge m}
		 |(h'(\tilda b_n)-h'(b))(\bar B_n-b)|\ge a\} \cr
	 &\le {K'|h'(b)|^{2+\lambda}\tau^{2+\lambda}
		\over a^{2+\lambda}}
	  +{\rm Pr}\{\sqrt{m}\sup_{n\ge m}\rho(|\bar B_n-b|)
			 |\bar B_n-b|\ge a\}, \cr}$$
where $K'$ is a new constant and 
writing $\rho(r)$ for the maximum of $|h'(x)-h'(b)|$ as 
$|x-b|\le r$. Let $\eps_0$ be such that $\rho(r)\le1$ when
$r\le\eps_0$ (we even have $\rho(r)\rightarrow0$ as $r\rightarrow0$),
and let $g(r)=\rho(r)r$, a continuously increasing function.
The second term above is bounded by 
$${\rm Pr}\{\sup_{n\ge m}|\bar B_n-b|\ge g^{-1}(a/\sqrt{m})\}
	 \le {K'\tau^{2+\lambda}
		\over \{\sqrt{m}g^{-1}(a/\sqrt{m})\}^{2+\lambda}},$$
which again is bounded by $K'\tau^{2+\lambda}/a^{2+\lambda}$, 
provided $\sqrt{m}g^{-1}(a/\sqrt{m})\ge a$, 
or $a/\sqrt{m}\ge g(a/\sqrt{m})$,
or $1\ge \rho(a/\sqrt{m})$. 
But this holds when $a/\sqrt{m}\le \eps_0$,
which proves (6.2). \firkant

\smallskip
This result extends without serious difficulties to $p$-dimensional
$\hatt\theta_n$ being a smooth function of $r$ averages.
With notation as in 3A the proviso for correct convergence 
of $E\eps^2N_\eps$ is finiteness of $E|B_{i,j}|^{2+\lambda}$ for
some positive $\lambda$, for $j=1,\ldots,r$. 
One may also look for conditions in the maximum likelihood 
estimator case. The essential
requirement is 
$E|\dell\log f(X_i,\theta_0)/\dell\theta_j|^{2+\lambda}<\infty$,
for $j=1,\ldots,p$. 

\bigskip
{\bf 7. The number of $\eps$-misses.}
We have been able to reach rather general and elegant results for
$N_\eps$ by the stochastic process approach, working with
$\sqrt{m}(\hatt\theta_{[mt]}-\theta_0)$ and its limiting 
process $\sigma_0W(t)/t$. 
This approach can also successfully be applied 
to other random non-observable quantities of interest,
thereby broadening the perspective.

\smallskip
{\sl 7A. The one-dimensional case.} 
To illustrate this point, 
consider $Q_\eps(a)$, the number of times, among $n\ge a/\eps^2$, 
where $|\hatt\theta_n-\theta_0|\ge\eps$. Then 
$$\eps^2Q_\eps(a)\rightarrow_d \sigma_0^2Q(a/\sigma_0^2)
	=\sigma_0^2\,\mu\{t\ge a/\sigma_0^2\colon|W(t)/t|\ge 1\}, \eqno(7.1)$$ 
in which $\mu$ is Lebesgue measure on the {half}{line}. 
This can be proved as follows, under conditions (2.2)--(2.3).
Write $Q_\eps(a)$ cleverly as
$\int_{\langle a/\eps^2\rangle}^\infty 
I\{|\hatt\theta_{[s]}-\theta_0|\ge\eps\}\,ds$,
then let $m=1/\eps^2$, and change to $t=s/m$. 
After tending to details similar to those of Section 2 the result is 
$$\eps^2Q_\eps(a)=\int_{\langle ma\rangle/m}^\infty
	I\{\sqrt{m}|\hatt\theta_{[mt]}-\theta_0|\ge1\}\,dt 
	\rightarrow_d \int_a^\infty I\{\sigma_0|W(t)/t|\ge 1\}\,dt,$$
and the limit can be rewritten as $\sigma_0^2Q(a/\sigma_0^2)$ above. 
%
%
There is also simultaneous convergence in distribution 
of $\bigl(\eps^2N_\eps,\eps^2Q_\eps(a)\bigr)$ to 
$\sigma_0^2\bigl(\sup_{t\ge1}|W(t)/t|^2,Q(a/\sigma_0^2)\bigr)$. 
This follows by measurability 
and a.s.~continuity of the appropriate functionals on
$D[a,b]$-spaces, and an extra argument to take care of 
the tail. It can also be proved via the continuous mapping
theorem on the function space on $[0,\infty)$ described in 
Remark (i) of Section 2. 
These results can also be proved for $a=0$. 
We leave the details out, but regularity conditions (2.2)--(2.3) 
suffice once more. 
In particular $\eps^2$ times the total number of $\eps$-misses 
goes to $\sigma_0^2Q(0)=\sigma_0^2\,\mu\{t\ge0\colon|W(t)/t|\ge1\}$.  
Note that $EQ(b)=E(\chi_1^2-b)I\{\chi_1^2\ge b\}$, 
using Fubini's theorem. In particular 
$EQ(0)=1$ and $EQ(0.95)=\half$, which means that 
the estimator sequence has about 
$\half/\eps^2$ misses of size $\sigma_0\eps$  
for $n\le 0.95/\eps^2$ and about $\half/\eps^2$ misses of size $\sigma_0\eps$
for $n\ge 0.95/\eps^2$. 
We mention finally that 
Kao (1978) has the $\eps^2Q_\eps(a)$ result, 
but only for the special case of simple i.i.d.~averages and $a=0$. 

\smallskip
{\sl 7B. The multi-dimensional case.} 
One result is the following, under conditions (2.6)--(2.7): 
Let $Q_\eps(a)$ be the 
number of times, among $n\ge a/\eps^2$, where  
$(\hatt\theta_n-\theta_0)'\Sigma_0^{-1}(\hatt\theta_n-\theta_0)\ge\eps^2$,
with notation and conditions as in the Theorem of Section 2.
Then $\eps^2Q_\eps(a)$ tends to 
$Q(a)=\mu\{t\ge a\colon \sum_{i=1}^pW_i(t)^2>t^2\}$. 
Note that $Q(a)$ has mean value $E(\chi^2_p-a)I\{\chi^2_p\ge a\}$,
which is easy to compute. In particular the total number 
of $\eps$-misses for the estimator sequence 
(with the Mahalanobis distance) is about $p/\eps^2$. 

\eject 

Another result with two interesting consequences is as follows:
Let distance function and conditions be as in the Theorem of Section 2, 
and let $Q_\eps$ be the total number of 
$\|\hatt\theta_n-\theta_0\|\ge\eps$ cases. 
Then $\eps^2Q_\eps$ tends to 
$Q=\mu\{t\ge0\colon\|\Sigma_0^{1/2}W(t)/t\|\ge1\}$.
Our first point is yet another asymptotic optimality 
property for the maximum likelihood sequence: 
In the limit, as $\eps\rightarrow0$,
provided the underlying parametric model is correct,  
{\sl no other estimator
sequence has stochastically fewer $\eps$-misses!}   
Our second point is that $EQ$ can be computed and leads to 
an a.r.e.~measure in the multi-dimensional case, 
cf.~the discussion that led to (2.5) and the end remark of Section 2.
Taking $\|x-y\|=\{(x-y)'A(x-y)\}^{1/2}$ the mean of $Q$ is 
$\int_0^\infty{\rm Pr}\{W(t)'\Sigma_0^{1/2}A\Sigma_0^{1/2}W(t)/t\ge t\}\,dt$,
which becomes $EZ'\Sigma_0^{1/2}A\Sigma_0^{1/2}Z$, where $Z\sim N_p(0,I_p)$.
Hence $EQ={\rm Tr}(A\Sigma_0)$. 
One can also prove convergence of $\eps^2EQ_\eps$ to $EQ$,
under conditions that in fact are simpler than those of Section 6. 
Suppose $\sqrt{n}(\hatt\theta_{n,j}-\theta_0)\rightarrow_dN(0,\Sigma_j)$
for $j=1,2$, and let $Q_{\eps,j}$ be the number of $\eps$-misses for  
method $j$. Then the arguments presented before (2.5) suggest 
$${\rm a.r.e.}=\lim_{\eps\rightarrow0}{EQ_{\eps,1}\over EQ_{\eps,2}}
	={{\rm Tr}(A\Sigma_1)\over {\rm Tr}(A\Sigma_2)}. \eqno(7.2)$$
Under ordinary Euclidian distance a.r.e.~becomes 
${\rm Tr}(\Sigma_1)/{\rm Tr}(\Sigma_2)$. See also 8E. 

\smallskip
{\sl 7C. The number of $\eps$-misses for Glivenko--Cantelli.}
Consider the more complicated situation of Section 4.
Let $Q_\eps$ be the number of times $\|F_n-F\|\ge\eps$. 
Combining arguments above with those of Section 4 one can show that 
$\eps^2Q_\eps\rightarrow_d Q=\mu\{s\colon A(s)\ge1\}$, in which
$A(s)=\max_{0\le t\le 1}|K_0(s,t)/s|$. But, for fixed $s$,
$K_0(s,.)/s$ is distributed like $W^0(.)/\sqrt{s}$, where $W^0(.)$
is a Brownian bridge, so that $A(s)=_d\|W^0\|/\sqrt{s}$, where
$\|W^0\|$ is the maximum of $|W^0(t)|$. This leads to 
$EQ=\int_0^\infty {\rm Pr}\{\|W^0\|\ge\sqrt{s}\}\,ds=E\|W^0\|^2=\pi^2/12$.
Accordingly the full estimator sequence will have about $0.822/\eps^2$ 
cases of $\|F_n-F\|\ge\eps$. --- Similarly, $\eps^2$ times the
total number of cases of $\int(F_n-F)^2\,dF\ge\eps^2$ will
converge in distribution to a variable with expected value $1/6$.        
And for a final example of a non-trivial result reached using
these methods, let $Q_\eps^*$ denote the number of $\int|F_n-F|\,dF\ge\eps$
cases. Then $\eps^2Q_\eps^*$ tends to an appropriate $Q^*$, 
and $\eps^2EQ^*_\eps$ tends to $EQ^*$, which can be proved to be 
equal to $E\bigl\{\int_0^1|W^0(t)|\,dt\bigr\}^2$, 
and which is found to be $7/60$ by stamenous calculations. 

\smallskip
{\sl 7D. The number of $\eps$-misses for a density estimator.}
Let finally $Q_\eps$ be the total number of times $|f_n(x)-f(x)|\ge\eps$,
in the density estimation problem considered in Section 5.
Analysis analogous to that above leads to 
$\eps^{5/2}Q_\eps\rightarrow_d Q=\mu\{s\colon |Z(s)|\ge1\}$,
where $Z(.)$ is the process defined in Section 5. 
One can then show that 
$$EQ=E\big|\half c^2f''(x)
  	+c^{-1/2}f(x)^{1/2}\beta_K^{1/2}N(0,1)\big|^{5/2}. \eqno(7.3)$$ 
The value of $c$ that gives 
the smallest expected number of $\eps$-misses,
in the limit as $\eps\rightarrow0$,  
can be shown to be $1.008\,c_0(x)$, as in Section 5.
Similar but more cumbersome calculations can be carried out 
for $Q'_\eps$, the number of times $|f_n'(x)-f'(x)|\ge\eps$, 
under the optimal scheme $h_n=cn^{-1/7}$. 
One finds that $\eps^{7/2}Q'_\eps$ tends to a certain $Q'$.
The best value of $c$ from the point of view of approximate
mean squared error is $c_0(x)=\{3\gamma_Kf(x)/f'''(x)^2\}^{1/7}$,
where $\gamma_K=\int K'(u)^2\,du$. But the value of $c$ that
minimises $EQ'$ can by determined efforts be shown to be $1.049\,c_0(x)$. 

\bigskip
{\bf 8. Complementary remarks and results.}

\smallskip
{\sl 8A. Numerical information.}
Central in our limit distribution results is the variable 
$W_{\max}=\max_{0\le s\le1}|W(s)|$. Its distribution 
can be found in Shorack and Wellner (1986, p.~35), for example. 
One can prove that 
$EW_{\max}=\sqrt{\pi/2}=1.2533$;
$EW_{\max}^2=2G=1.8319$, featuring Catalan's constant; 
${\rm Var}\,W_{\max}=2G-\half\pi=0.5110^2$; 
${\rm stdev}(W_{\max}^2)=(EW_{\max}^4-4G^2)^{1/2}=1.6055$. 
In the case of a single parameter, therefore, the following holds,
in the notation of Section 6:
$\eps\,EN_\eps^{1/2}\rightarrow\sqrt{\pi/2}\,\sigma_0$,
if only $E|Z_i|^2<\infty$;
$\eps^2\,EN_\eps\rightarrow2G\sigma_0^2$, 
if $E|Z_i|^{2+\lambda}$ is finite for some positive $\lambda$;
$\eps^2\,{\rm stdev}(N_\eps)\rightarrow 1.6055\,\sigma_0^2$,
if $E|Z_i|^{4+\lambda}$ is finite. The distribution of $N_\eps$ 
is skewed to the right, as 
${\rm skew}(N_\eps)=E\{(N_\eps-EN_\eps)/{\rm stdev}(N_\eps)\}^3
\rightarrow2.3308$ if $E|Z_i|^{6+\lambda}$ is finite. 
%
In the case of several parameters and the Mahalanobis
distance we have proved 
$\eps^2N_\eps\rightarrow_d \chi^2_{p,\max}$,
the maximum of $\chi^2_p(s)=\sum_{i=1}^pW_i(s)^2$ over $[0,1]$.
A way of computing its distribution is provided by DeLong (1980),
along with a few quantiles. 
[More details and a fuller table of quantiles,
arrived at by simulation of Brownian motions, 
are given in Hjort and Fenstad (1990).]

\smallskip 
{\sl 8B. Extension to non-i.i.d.~situations.} 
Our basic results read 
$\eps^2N_\eps\rightarrow_d\sigma_0^2W_{\max}^2$
and $\eps^2Q_\eps(0)\rightarrow\sigma_0^2Q(0)$
(in the one-parameter case), where $\sigma_0^2$ is the variance of the
limit distribution for $\sqrt{m}(\hatt\theta_m-\theta_0)$.
These continue to hold for large classes non-i.i.d.-situations. 
The key ingredients are process convergence 
$\sqrt{m}(\hatt\theta_{[mt]}-\theta_0)\rightarrow_d \sigma_0W(t)/t$
in $D[b,c]$-spaces 
(tightness and convergence of finite-dimensional distributions)
and a tail inequality for $[c,\infty)$. 
Proving this for a particular case requires 
attention to technical details depending upon that case, however.  
In the technical report version of this article 
such attention is given to linear regression and 
to a situation with auto-correlation.

\smallskip
{\sl 8C. A slow minimax estimator.}
Let $X_1,X_2,\ldots$ be independent Bernoulli trials with 
success probability $p$. The maximum likelihood 
estimator for $p$ after $n$ trials is $\hatt p_n=Y_n/n$, 
where $Y_n$ is the number or successes in the first $n$ trials.
From earlier results we know that
$\eps^2N_\eps\rightarrow_d p(1-p)W_{\max}^2$,
where $N_\eps$ is the last time $|\hatt p_n-p|\ge\eps$. --- Now
consider the minimax estimator 
$p^*_n=(\sqrt{n}\hatt p_n+\half)/(\sqrt{n}+1)$, and the accompanying
$N_\eps^*$, the last time $|p^*_n-p|\ge\eps$. 
Some analysis reveals that  
$$\sqrt{m}(p^*_{[mt]}-p) 
	 \rightarrow_d\sqrt{p(1-p)}\Bigl[{W(t)\over t}
		+{\half-p\over \sqrt{p(1-p)}}{1\over \sqrt{t}}\Bigr]
		{\rm \ in\ } D[1,c]. $$
This can be used to prove $\eps^2N_\eps^*\rightarrow_d p(1-p)
\max_{0\le s\le1}|W(s)+b(p)\sqrt{s}|^2$, where 
$b(p)=(\half-p)/\{p(1-p)\}^{1/2}$.
Accordingly $N_\eps^*$ for $p^*_n$ is stochastically larger, in the limit,
than $N_\eps$ for $\hatt p_n$ (unless $p={1\over2}$). 
%
There is a similar story for $Q_\eps$ and $Q_\eps^*$, 
the number of times $\hatt p_n$ and $p_n^*$ miss with more than $\eps$.
One can prove that $\eps^2Q_\eps\rightarrow_dQ$
and $\eps^2Q^*_\eps\rightarrow_dQ^*$, where 
$EQ=p(1-p)$ and $EQ^*=p(1-p)+(\half-p)^2={1\over4}$. 

There are analogous results for the cdf-estimator 
$F_n^*=(\sqrt{n}F_n+\half)/(\sqrt{n}+1)$, which can be shown 
to be minimax when the loss function is $\int(\hatt F-F)^2\,dw$, 
$w$ any given weight function with mass 1. 
Then $F_n^*$ can expect ${1\over4}/\eps^2$ instances with 
loss $\ge\eps^2$, regardless of the underlying $F$, 
whereas the non-minimax estimator $F_n$ 
can expect $\int F(1-F)\,dw/\eps^2$ such instances. 

\smallskip
{\sl 8D. Other distances.} 
Our basic result (2.8) was phrased in terms of a distance 
function $\|\hatt\theta_n-\theta_0\|$. 
The arguments carry through also for other measures of distance
that are not of the norm type. As a particular example of some
interest, let 
$d[\theta_0\colon\theta]$
be the Kullback--Leibler distance 
$\int f(x,\theta_0)\log\{f(x,\theta_0)/f(x,\theta)\}\,dx$, 
in some model with a $p$-dimensional parameter. 
Let $\hatt\theta_n$ be the maximum likelihood
estimator and let $M_\eps$ be the last $n$ at which
$d[\theta_0\colon\hatt\theta_n]\ge\eps$.
Then $2\eps M_\eps\rightarrow_d \chi^2_{p,\max}$ of (2.10)
can be proved under mild conditions. 
Note that the limit is the same regardless of the actual 
parametric family.
This holds when $f(x,\theta_0)$ represents the true model.
A more general result, 
valid under the agnostic viewpoint presented in Section 3B,
is given in Hjort and Fenstad (1990), along with 
further examples with other distance functions.


\smallskip
{\sl 8E. Second order results.} 
Our a.r.e.~measures in (2.5) and (7.2) do not distinguish
between estimators with the same limiting distribution.
To do so requires second order asymptotics for 
$\eps^2N_\eps$ and $\eps^2Q_\eps$. 
In Hjort and Fenstad (1991) and Hjort (1991) we have sorted out 
the limiting behaviour of differences between $Q_\eps$'s in cases 
where their ratio tends to 1, thereby making it possible to 
find second order optimal estimator sequences in many cases of interest.  
Thus, in the binomial situation, 
the $(Y_n+{2\over3})/(n+{4\over3})$ sequence can be expected 
to make 2.667 fewer $\eps$-errors than the traditional $Y_n/n$ sequence,
for example, regardless of the underlying $p$ parameter.
And among all estimators of the form   
$\sum_{i=1}^n(X_i-\bar X_n)^2/(n+c)$ for a normal variance 
the one with denominator $n-{1\over3}$ 
can be expected to make the fewest $\eps$-errors.


\smallskip 
{\sl 8F. Sequential fixed-volume confidence regions.}
Suppose (2.5) and (2.6) hold, and write 
$N^*_\eps$ for the last time 
$(\hatt\theta_n-\theta_0)'\hatt\Sigma_n^{-1}
(\hatt\theta_n-\theta_0)\ge\eps^2$. 
Then $\eps^2N^*_\eps$ tends to $\chi^2_{p,\max}$ of (2.10),
provided merely that 
$\hatt\Sigma_n\rightarrow\Sigma_0$ a.s.~(convergence 
in probability does not suffice). 
Let $\eps$ be small and given, 
find $c$ such that ${\rm Pr}\{\chi^2_{p,\max}\le c\}=0.95$, 
put $m=[c/\eps^2]$, and consider  
$I_n^*=\{\theta\colon(\theta-\hatt\theta)'\hatt\Sigma_n^{-1}
	  (\theta-\hatt\theta_n)\le\eps^2\}$. 
Then 
${\rm Pr}\{\theta_0\in I_n^*{\rm \ for \ all\ }n\ge m\}\doteq0.95$. 
The details of this construction are in Hjort and Fenstad (1990). 

\smallskip
{\sl 8G. Shrinking boundaries and tests with power 1.} 
Methods of this paper can be used to construct 
sequential confidence regions with shrinking volume 
as well as sequential tests with power 1. 
See Hjort and Fenstad (1990). 

\smallskip
{\sl 8H. The probabilities of $N_{\eps,1}<N_{\eps,2}$
and $Q_{\eps,1}<Q_{\eps,2}$.} 
Consider two competing estimator sequences, 
with accompanying last $\eps$-miss variables $N_{\eps,j}$ 
and number of $\eps$-misses variables $Q_{\eps,j}$, 
as in (2.5) and (7.2). 
The probabilities ${\rm Pr}\{N_{\eps,1}<N_{\eps,2}\}$
and ${\rm Pr}\{Q_{\eps,1}<Q_{\eps,2}\}$ 
will usually converge as $\eps$ goes to zero;
in fact $\eps^2(N_{\eps,1},N_{\eps,2},Q_{\eps,1},Q_{\eps,2})$
has a joint limiting distribution in terms of two 
correlated Brownian motions, under natural conditions. 
These limits are found in Hjort and Fenstad (1990). 
As an example, consider the average estimator 
$\hatt\theta_{n,1}=\bar X_n$ and 
the median estimator $\hatt\theta_{n,2}=M_n$ 
for the mean parameter in the normal model. 
Then $N_{\eps,1}<N_{\eps,2}$ with probability about 0.72,
and $Q_{\eps,1}<Q_{\eps,2}$ with probability about 0.69. 
To give $\bar X_n$ a harder match, replace the second estimator
with the solution of $\sum_{i=1}^n\arctan(X_i-\theta)=0$, 
an $M$-estimator with a smooth and bounded influence function. 
Then the figures become respectively 0.56 and 0.55.
These are simulation-based figures computed using 
the exact limit distributions. 

\bigskip
\baselineskip11pt
{\bf Acknowledgements.} 
We are grateful to a reviewer and to David Siegmund 
for comments that inspired improvements and for pointing out 
some of the probabilistic literature which we had missed 
in an earlier version. 

\vfill \eject 

\bigskip
\centerline{\bf References}

\medskip\noindent\parskip2pt\parindent0pt\baselineskip11pt
\rm 

Bahadur, R.R. (1967).
Rates of convergence of estimates and test statistics.
{\sl Ann.~Math. Statist.}~{\bf 38}, 303--324.

von Bahr, B. (1965). 
On convergence of moments in the central limit theorem.
{\sl Ann. Math.~Statist.} {\bf 36}, 808--818.


Beran, R. (1977). 
Estimating a distribution function.
{\sl Ann.~Statist.}~{\bf 5}, 400--404. 
 
Bickel, P.J.~and Wichura, M.J. (1971). 
Convergence criteria for multiparameter stochastic
processes and some applications.
{\sl Ann.~Math.~Statist.}~{\bf 42}, 1656--1670.

Bickel, P.J.~and Rosenblatt, M. (1973). 
On some global measures of the deviations of density function 
estimates.
{\sl Ann.~Statist.}~{\bf 1}, 1071--1095. 
Corrigenda {\sl ibid.} (1975, p.~1370).

Billingsley, P. (1968).
{\sl Convergence of Probability Measures.}
Wiley, New York.



DeLong, D.M. (1980).
Some asymptotic properties of a progressively 
censored nonparametric test for multiple regression.
{\sl J.~Mult.~Anal.}~{\bf 10}, 363--370.





Hjort, N.L. (1986). 
{\sl Statistical Symbol Recognition.}
Research monograph, Norwegian Computing Centre, Oslo. 

Hjort, N.L. (1991). 
Some exponential distributions associated with Brownian motion.
Statistical Research Report, Dept.~of Math., University of Oslo, 
submitted for publication. 

Hjort, N.L.~and Fenstad, G. (1990).
On the last time a strongly consistent estimator
is more than $\eps$ from its target value.
Statistical Research Report, Dept.~of Math., University of Oslo.

Hjort, N.L.~and Fenstad, G. (1991).
Some second order asymptotics for the number of times 
an estimator is more than $\eps$ from its target point.
Statistical Research Report, Dept.~of Math., University of Oslo;
submitted for publication. 


Kao, C-S. (1978). 
On the time and the excess of linear boundary crossings of sample sums.
{\sl Ann.~Statist.}~{\bf 6}, 191--199.


Lehmann, E.L. (1983).
{\sl Theory of Point Estimation.}
Wiley, New York. 

Lo\`eve, M. (1960).
{\sl Probability Theory.}
2nd Ed., Van Nostrand, Toronto. 

M\"uller, D.W. (1968).
Verteilungs-Invarianzprinzipien f\"ur das starke Gesetz 
der gro\ss en Zahl. 
{\sl Z.~Wahrscheinlichkeitstheorie verw.~Geb.}~{\bf 10}, 173--192.

M\"uller, D.W. (1970). 
On Glivenko--Cantelli convergence.
{\sl Z.~Wahrscheinlichkeitstheorie verw.~Geb.}~{\bf 16}, 195--210.

M\"uller, D.W. (1972).
Randomness and extrapolation.
Proceedings of the Sixth Berkeley Symposium, II, 1--31.



Robbins, H., Siegmund, D., and Wendel, J. (1968). 
The limiting distribution of the last time $s_n\ge n\eps$.
{\sl Proc.~Nat.~Acad.~Sci.~USA}~{\bf 61}, 1228--1230.


Serfling, R. (1980). 
{\sl Approximation Theorems of Mathematical Statistics.}
Wiley, New York. 

Shao, J. (1989).
Functional calculus and asymptotic theory for statistical analysis. 
{\sl Statist. and Probab.~Letters} {\bf 8}, 397--405. 


Shorack, G.R.~and Wellner, J.A. (1986).
{\sl Empirical processes with applications to statistics.}
Wiley, Singapore. 

Stute, W. (1983). 
Last passage time of $M$-estimators.
{\sl Scand.~J.~Statist.}~{\bf 10}, 301--305.



\bye